\documentclass[10pt]{article}      

\usepackage{amsmath, amssymb, amsthm, amscd}   

\usepackage{graphicx}                   
\usepackage{psfrag}                     

\input{kleinmacros.sty}
\input{kleinthm.sty}

\input{kleindefs.sty}

\newcommand{\matR} {\ensuremath {\mathbb{R}}}
\newcommand{\matQ} {\ensuremath {\mathbb{Q}}}
\newcommand{\matH} {\ensuremath {\mathbb{H}}}
\newcommand{\matZ} {\ensuremath {\mathbb{Z}}}

\newcommand{\matP} {\ensuremath {\mathbb{P}}}
\newcommand{\matRP} {\ensuremath {\mathbb{RP}}}
\newcommand{\calM} {\ensuremath {\mathcal{M}}}
\newcommand{\calH} {\ensuremath {\mathcal{H}}}

\newcommand{\Sol}{{\rm Sol}}

\newcommand{\Vol}{{\rm Vol}}
\newcommand{\timtil}{\begin{picture}(12,12)
\put(2,0){$\times$}\put(2,4.5){$\sim$}\end{picture}}

\newcommand{\chiorb}{\chi^{\rm orb}}
\newcommand{\voltet}{v_{\rm T}}
\newcommand{\voloct}{v_{\rm O}}
\newcommand{\voldrum}{v_{\rm D}}
\newcommand{\ptwoirred}{$\matP^2$-irreducible}

\newcommand{\JSJ}{{\rm JSJ}}
\newcommand{\matr} [4] 
{\left(\tiny{\begin{array}{@{}c@{\ }c@{}} #1 & #2 \\ #3 & #4 \\ \end{array}} \right)}
\newcommand{\tr} {{\rm tr\,}}

\begin{document}

\renewcommand{\title}{Complexity of $3$-manifolds}
\renewcommand{\author}{Bruno Martelli\footnote{Supported by the INTAS project ``CalcoMet-GT'' 
03-51-3663}}

\dohead

\markboth{Martelli}{Complexity}   

\renewcommand{\firstpage}{\pageref{martelli_first}}  
\renewcommand{\lastpage}{\pageref{martelli_last}}

\begin{abstract}
We give a summary of known results on Matveev's complexity of compact $3$-manifolds.
The only relevant new result is
the classification of all closed orientable irreducible $3$-manifolds
of complexity $10$.
\end{abstract}

\section{Introduction} \label{martelli_first}  
In 3-dimensional topology, various quantities are defined, that measure how complicated a 
compact 3-manifold $M$ is. 
Among them, we find the Heegaard genus, the minimum number of tetrahedra in a triangulation,
and Gromov's norm (which equals the volume when $M$ is hyperbolic).
Both Heegaard genus and Gromov norm are additive on connected sums, and behave well
with respect to other common cut-and-paste operations,
but it is hard to classify all manifolds with a given genus or norm.
On the other hand, triangulations with $n$ tetrahedra are more suitable for computational
purposes, since they are finite in number
and can be easily listed using a computer, 
but the minimum number of tetrahedra is
a quantity which does not behave well with any cut-and-paste operation on 3-manifolds.
(Moreover, it is not clear what is meant by ``triangulation'':
do the tetrahedra need to be embedded? Are ideal vertices admitted when $M$ has boundary?)

In 1988, Matveev introduced~\cite{Mat88} for any compact $3$-manifold $M$ a non-negative integer
$c(M)$, which he called the \emph{complexity} of $M$, defined as the minimum number
of vertices of a \emph{simple spine} of $M$. The function $c$ is finite-to-one on the
most interesting sets of compact 3-manifolds, and it behaves well with respect to the
most important cut-and-paste operations. Its main properties are listed below.
\begin{description}
\item[additivity] $c(M\#M')=c(M)+c(M')$;
\end{description}
\begin{description}
\item[finiteness] for any $n$ there is a finite number of closed \ptwoirred\ $M$'s
with $c(M)=n$, and a finite number of hyperbolic $N$'s with $c(N) = n$;
\end{description}
\begin{description}
\item[monotonicity] $c(M_F)\leqslant c(M)$ for any incompressible $F\subset M$ cutting $M$ into $M_F$.
\end{description}

We recall some definitions used throughout the paper. 
Let $M$ be a compact $3$-manifold, possibly with boundary.
We say that $M$ is \emph{hyperbolic} if it admits (after removing all tori 
and Klein bottles from the boundary) 
a complete hyperbolic metric of finite volume (possibly with cusps and geodesic
boundary). Such a metric is unique by Mostow's theorem (see~\cite{McMu} for a proof).
A surface in $M$ is \emph{essential} if it is incompressible, $\partial$-incompressible,
and not $\partial$-parallel.
Thurston's Hyperbolicity 
Theorem for Haken manifolds ensures that a compact $M$ with boundary 
is hyperbolic if and only if
every component of $\partial M$ has $\chi\leqslant 0$, and $M$
does not contain essential surfaces with $\chi\geqslant 0$.
The complexity satisfies also the following strict inequalities.
\begin{description}
\item[filling] every closed hyperbolic $M$ 
is a Dehn filling of some hyperbolic $N$ with $c(N)<c(M)$;
\end{description}
\begin{description}
\item[strict monotonicity] $c(M_F)<c(M)$ if $F$ is essential and $M$ is closed
\ptwoirred\ or hyperbolic;
\end{description}

Some results in complexity zero already show that
the finiteness property does not hold for all compact $3$-manifolds.
\begin{description}
\item[complexity zero] the closed \ptwoirred\ manifolds with $c=0$ are $S^3, 
\matRP^3,$ and $L(3,1)$. We also have $c(S^2\times S^1) = c(S^2\timtil S^1) = 0$.
Interval bundles over surfaces and handlebodies also have $c=0$.
\end{description}

The ball and the solid torus have therefore complexity zero.
Moreover, the additivity property actually also holds for $\partial$-connected sums. 
These two facts together imply the following.
\begin{description}
\item[stability]
The complexity of $M$ does not change when adding 1-handles to $M$ or removing interior
balls from it.
\end{description}

Note that both such operations that not affect $c$ are ``invertible'' and hence
topologically inessential.
In what follows, a simplicial face-pairing $T$ of some tetrahedra is a 
\emph{triangulation} of a closed 3-manifold $M$ when $M = |T|$. Tetrahedra are therefore 
not necessarily embedded in $M$. A simplicial pairing $T$ is
an \emph{ideal triangulation} of a compact $M$ with boundary if $M$ is $|T|$ minus open stars of
all the vertices. The finiteness property above follows easily from the following.
\begin{description}
\item[naturality] if $M$ is closed \ptwoirred\ and not $S^3, \matRP^3,$ or $L(3,1)$,
then $c(M)$ is the minimum number of tetrahedra in a triangulation of $M$. If $N$ is
hyperbolic with boundary, 
then $c(N)$ is the minimum number of tetrahedra in an ideal triangulation of $N$.
\end{description} 

The beauty of Matveev's complexity theory relies on the fact that simple spines are more
flexible than triangulations: for instance spines can often be simplified by puncturing
faces, and can always be cut along normal surfaces. In particular, we have the following
result. An (ideal) triangulation $T$ of $M$ is \emph{minimal} when $M$ cannot be (ideally) 
triangulated with fewer tetrahedra. A \emph{normal surface} in $T$ is one intersecting
the tetrahedra in normal triangles and squares, see~\cite{normal}.
\begin{description}
\item[normal surfaces] let $T$ be a minimal (ideal) triangulation of a closed \ptwoirred\
(hyperbolic with boundary) manifold $M$ different from $S^3$, $\matRP^3$, 
and $L(3,1)$. If $F$ is a normal
surface in $T$ containing some squares, then $c(M_F)< c(M)$.
\end{description}

As an application of the previous properties, the following result was implicit in
Matveev's paper~\cite{Mat90}.
\begin{corollary}
Let $T$ be a minimal triangulation of a closed \ptwoirred\
3-manifold $M$ different from $S^3,\matRP^3, L(3,1)$. 
Then $T$ has one vertex only, and it contains no normal spheres, 
except the vertex-linking one.
\end{corollary}
Computers can easily handle spines and triangulations, 
and manifolds of low complexity have been classified
by various authors. 
Closed orientable irreducible 
manifolds with $c\leqslant 6$ were classified by Matveev~\cite{Mat88} 
in 1988. Those with $c=7$ were then classified in 1997 by 
Ovchinnikov~\cite{Ov, Mat:book}, 
and those with $c=8,9$ in 2001 by Martelli and Petronio~\cite{MaPe}.
We present here the results we recently found for $c=10$. 
The list of all manifolds with $c=10$ has
also been computed independently by Matveev~\cite{Mat:priv}, and the two
tables (each consisting of $3078$ manifolds) coincide.
The closed \ptwoirred\ non-orientable manifolds with $c\leqslant 7$ have been
listed independently in 2003 by
Amendola and Martelli~\cite{AmMa2}, and Burton~\cite{Bu}.

Hyperbolic manifolds with cusps and without geodesic boundary were listed for 
all $c\leqslant 3$ in the orientable case by Matveev and Fomenko~\cite{MaFo} in 1988, and for 
all $c\leqslant 7$ by Callahan, Hildebrand, and Weeks~\cite{CaHiWe} in 1999. Orientable hyperbolic
manifolds with geodesic boundary (and possibly some cusps) were listed for
$c\leqslant 2$ by Fujii~\cite{Fuj} in 1990, and for
$c\leqslant 4$ by Frigerio, Martelli, and Petronio~\cite{FriMaPe2} in 2002.

All properties listed above were proved by Matveev in~\cite{Mat90}, and extended
when necessary to the non-orientable case by Martelli and Petronio 
in~\cite{MaPe:nonori}, except the filling property, which is a new result proved below
in Subsection~\ref{properties:subsection}.
The only other new results contained in this paper are the complexity-$10$ closed census
(also constructed independently by Matveev~\cite{Mat:priv}),
and the following counterexample (derived from that census) of a conjecture
of Matveev and Fomenko~\cite{MaFo} stated in Subsection~\ref{counterexample:subsection}.
\begin{proposition}
There are two closed hyperbolic fillings $M$ and $M'$ of the same cusped hyperbolic $N$
with $c(M)<c(M')$ and $\Vol(M)>\Vol(M')$.
\end{proposition}
We mention the most important discovery of our census.
\begin{proposition}
There are $25$ closed hyperbolic manifolds with $c=10$ (while none with $c\leqslant 8$
and four with $c=9$).
\end{proposition}

This paper is structured as follows: the complexity of a $3$-manifold is defined in 
Section~\ref{definitions:section}. We then collect
in Section~\ref{closed:section} and~\ref{hyperbolic:section} the censuses of
closed and hyperbolic 3-manifolds described above, together with the new results in
complexity $10$. Relations between
complexity and volume of hyperbolic manifolds are studied in 
Section~\ref{complexity:volume:section}.
Lower bounds for the complexity, together with some infinite families of 
hyperbolic manifolds with boundary for
which the complexity is known, 
are described in Section~\ref{lower:bounds:section}. 
The algorithm and tools usually employed to produce a census are described in 
Section~\ref{minimal:section}.
Finally, we describe the decomposition of a manifold into \emph{bricks} 
introduced by Martelli and Petronio in~\cite{MaPe, MaPe:nonori}, necessary for our closed census
with $c=10$, in Section~\ref{bricks:section}. All sections may be read
independently, except that Sections~\ref{minimal:section} and~\ref{bricks:section}
need the definitions contained in Section~\ref{definitions:section}.

\subsection*{Acknowledgements} 
The author warmly thanks Carlo Petronio and Sergej Matveev for their continuous support.
He also thanks the referee for her/his suggestions.

\section{The complexity of a 3-manifold} \label{definitions:section}
We define here simple and special spines, and the complexity of a 3-manifold.
We then show a nice relation between spines without vertices
and Riemannian geometry, found by 
Alexander and Bishop~\cite{AleBi}. Finally, we prove the filling property stated in the
Introduction.

\subsection{Definitions} \label{definitions:subsection}
\begin{figure}
\centrefig{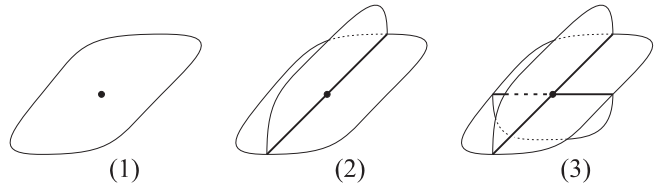}
\caption{Neighborhoods of points in a special polyhedron.} 
\label{special:fig}
\end{figure}

We start with the following definition.
\begin{defn}
A compact 2-dimensional
polyhedron $P$ is \emph{simple} if the link of every point in $P$ is contained in
the graph \includegraphics[width = .4cm]{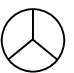}.
\end{defn}
Alternatively, $P$ is simple if it is locally contained in the 
polyhedron shown in Fig.~\ref{special:fig}-(3).
A point, a compact graph, a
compact surface are therefore simple. 
The polyhedron given by two orthogonal discs intersecting
in their diameter is not simple.
Three important possible kinds of
neighborhoods of points are shown in Fig.~\ref{special:fig}. 
A point having the whole of \includegraphics[width =.3 cm]{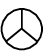} as a link is 
called a \emph{vertex}, and its regular neighborhood
is shown in Fig.~\ref{special:fig}-(3). 
The set $V(P)$ of the vertices of $P$ consists of isolated points, so
it is finite. Note that points, graphs, and surfaces do not contain vertices.

\begin{figure}
\centrefig[.6]{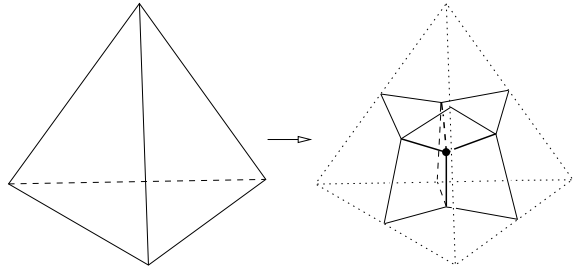}
\caption{A special spine of $M$ is dual to a triangulation, which is ideal or 1-vertex, depending
on whether $M$ has boundary or not.}
\label{dualspine:fig}
\end{figure}

A compact polyhedron $P\subset M$ is a \emph{spine} of a compact manifold $M$ with boundary
if $M$ collapses onto $P$. When $M$ is closed, we say that $P\subset M$ is a \emph{spine}
if $M\setminus P$ is an open ball. 
\begin{defn}
The \emph{complexity} $c(M)$ of
a compact 3-manifold $M$ is the minimal number of vertices of
a simple spine of $M$.
\end{defn}
As an example, a point is a spine of $S^3$, and therefore $c(S^3)=0$.
A simple polyhedron is \emph{special} when
every point has a neighborhood of one of the types (1)-(3) shown in
Fig.~\ref{special:fig}, and the sets of such points induce a
cellularization of $P$. That is, defining $S(P)$ as the set of points of type (2) or (3),
the components of $P \setminus S(P)$ 
should be open discs -- the \emph{faces} -- and the components of $S(P)\setminus V(P)$ 
should be open
segments -- the \emph{edges}.
\begin{remark} \label{spine:triangulation:rem}
A special spine of a compact $M$ with boundary is dual to an ideal triangulation of $M$,
and a special spine of a closed $M$ is dual to a 1-vertex triangulation of $M$, as suggested
by Fig.~\ref{dualspine:fig}. In particular, a special spine is a spine of a unique
manifold. Therefore the naturality property of $c$ may be read as follows:
every closed irreducible or hyperbolic manifold distinct from $S^3,\matRP^3$, and $L(3,1)$
has a special spine with $c(M)$ vertices. Such a special spine is then called \emph{minimal}.
\end{remark}

\subsection{Complexity zero}
A handlebody $M$ collapses onto a graph, which has no vertices, hence $c(M)=0$.
An interval bundle $M$ over a surface has that surface as a spine, and hence $c(M)=0$ again.
Note that, by shrinking the fibers of the bundle,
the manifold $M$ admits product metrics with arbitrarily small
injectivity radius and uniformly bounded curvature. This is a particular case of
a relation between spines and Riemannian geometry found
by Alexander and Bishop~\cite{AleBi}.
A Riemannian 3-manifold $M$ is \emph{thin} when
its curvature-normalized injectivity radius is less than some constant $a_2\approx 0.075$,
see~\cite{AleBi} for details. We have the following.
\begin{proposition}[Alexander-Bishop~\cite{AleBi}]
A thin Riemannian 3-manifold has complexity zero.
\end{proposition}

\subsection{The filling property} \label{properties:subsection}
We prove here the filling property, stated in the Introduction. 
Recall from~\cite{Mat90, Mat:book} 
that by thickening a special spine $P$ of $M$ we get
a handle decomposition $\xi_P$ of the same $M$. 
Normal surfaces in $\xi_P$ correspond to normal 
surfaces in the (possibly ideal) triangulation dual to $P$.

\begin{theorem} \label{filling:teo}
Every closed hyperbolic manifold $M$ is a Dehn filling of some hyperbolic $N$ with
$c(N)<c(M)$.
\end{theorem}
\begin{proof}
Let $P$ be a minimal special spine of $M$, which exists by Remark~\ref{spine:triangulation:rem}.
Take a face $f$ of $P$. By puncturing $f$ and collapsing the resulting polyhedron
as much as possible, we get a simple spine $Q$ of some $N$ obtained by drilling $M$
along a curve. Since $P$ is special, $f$ is incident to at least one vertex.
During the collapse, 
all vertices adjacent to $f$ have disappeared, hence $Q$ has less vertices than $P$.
This gives $c(N)<c(M)$. 

If $N$ is hyperbolic we are done. Suppose it is not. Then it is 
reducible, Seifert, or toroidal.
If $N$ is reducible, the drilled
solid torus is contained in a ball of $M$ and we get $N=M\# M'$
for some $M'$, hence $c(M)\leqslant c(N)<c(M)$ by the additivity property.
Then $N$ is irreducible.
Moreover $\partial N$ is incompressible
in $N$ (because $M$ is not a lens space). 
Then the 1-dimensional portion of $Q$ can be removed, and we can suppose $Q\subset P$
is a spine of $N$ having only points of the type of Fig.~\ref{special:fig}.

Our $N$ cannot be Seifert (because $M$ is hyperbolic), hence its \JSJ\ decompostion
consists of some tori $T_1,\ldots, T_k$. Each $T_i$ is essential in $N$
and compressible in $M$. Each $T_i$ can be isotoped in normal position with respect
to $\xi_Q$. Since $Q\subset P$, every normal surface in $\xi_Q$ is normal
also in $\xi_P$. The only normal surface in $\xi_P$ not containing squares is the
vertex-linking sphere,
therefore we have $c(M_{T_i})< c(M)$ for all $i$ by the normal surfaces property.
Each $T_i$ is compressible in $M$, hence either it 
bounds a solid torus or is contained in a ball. The latter case is excluded, otherwise
$M_{T_i}$ is the union of $M\# M'$ and a solid torus, and $c(M)\leqslant c(M_{T_i})<c(M)$.

Therefore each $T_i$ bounds a solid torus in $M$. Each solid torus contains the
drilled curve, hence they all intersect, and there is a solid torus $H$ bounded by a $T_i$ 
containing all the others. Therefore $M_{T_i} = N' \cup H$ where
$N'$ is a block of the \JSJ\ decomposition, 
which cannot be Seifert, hence it is hyperbolic. We have $c(N')=c(M_{T_i})<c(M)$, 
and $M$ is obtained by filling $N'$, as required.
\end{proof}

\begin{remark} The proof Theorem~\ref{filling:teo} is also 
valid for $M$ \emph{hyperbolike}, \emph{i.e.}~irreducible, atoroidal, and not Seifert.
\end{remark}

\section{Closed census} \label{closed:section}
We describe here the closed orientable irreducible manifolds with $c\leqslant 10$, 
and the closed non-orientable \ptwoirred\ ones with $c\leqslant 7$.
Such manifolds are collected in terms of their geometry, if any, in Table~\ref{closed:table}. 
The complete list of manifolds can be downloaded from~\cite{weblist}.
\begin{table}
  \begin{center}
    \begin{tabular}{rccccccccccc}
      & $0$ & $1$ & $2$ & $3$ & $4$ & $5$ & $6$ & $7$ & $8$ & $9$ & $10$ \\
      
      \multicolumn{12}{c}{orientable\phantom{\Big|}} \\

      lens spaces &
      $3$ & $2$ & $3$ & $6$ & $10$ & $20$ & $36$ & $72$ & $136$ & $272$ & $528$ \\

      other elliptic &
      . & . & $1$ & $1$ & $4$ & $11$ & $25$ & $45$ & $78$ & $142$ & $270$\\

      flat &
      . & . & . & . & . & . & $6$ & .    & .    & .      & .	 \\
      
      Nil &
      . & . & . & . & . & . & $7$ & $10$ & $14$ & $15$   & $15$\\

      ${\rm SL}_2\matR$ &
      . & . & . & . & . & . & .   & $39$ & $162$ & $513$ & $1416$\\
      
      \Sol &
      . & . & . & . & . & . & .   & $5$  & $9$   & $23$  & $39$\\

      $\matH^2\times\matR$ &
      . & . & . & . & . & . & .   & .    & $2$   & .     & $8$\\
      
      hyperbolic &
      . & . & . & . & . & . & .   & .    & .     & $4$   & $25$\\
      
      not geometric &
      . & . & . & . & . & . & .   & $4$  & $35$  & $185$ & $777$\\
      
      total orientable  &
      $\bf{3}$ & $\bf{2}$ & $\bf{4}$ & $\bf{7}$ & $\bf{14}$ & $\bf{31}$ & $\bf{74}$ & $\bf{175}$ & $\bf{436}$ & $\bf{1154}$ & $\bf{3078}$ \\

      \multicolumn{12}{c}{non-orientable\phantom{\Big|}} \\

      flat &
      . & . & . & . & . & . & $4$     & .   & & & \\

      $\matH^2\times\matR$ &
      . & . & . & . & . & . & .       & $2$ & & & \\

      \Sol &
      . & . & . & . & . & . & $1$     & $1$ & & & \\

      total non-orientable &
      . & . & . & . & . & . & $\bf 5$ & $\bf 3$ & & & \\

    \end{tabular}
  \end{center} 
  \caption{The number of closed \ptwoirred\ manifolds of
    given complexity (up to $10$ in the orientable case, 
and up to $7$ in the non-orientable
one) and geometry. Recall that there is no \ptwoirred\ manifold of type $S^2\times\matR$, and
no non-orientable one of type $S^3$,  Nil, and ${\rm SL}_2\matR$.
}
  \label{closed:table}
\end{table}

\subsection{The first $7$ geometries}
We recall~\cite{Sco} that there are eight important 3-dimensional geometries, six of them 
concerning Seifert manifolds. 
A Seifert fibration is described via its \emph{normalized parameters}
$\big(F, (p_1, q_1), \ldots, (p_k, q_k), t\big)$, where $F$ is a closed surface, 
$p_i>q_i>0$ for all $i$, and $t\geqslant -k/2$
(obtained by reversing orientation if necessary).
The Euler characteristic $\chi^{\rm orb}$ of the base orbifold and the Euler number
$e$ of the fibration are given respectively by
$$\chi^{\rm orb} = \chi(F)-\sum_{i=1}^k \left(1-\frac 1{p_i}\right), 
\qquad e = t + \sum_{i=1}^k \frac{q_i}{p_i}$$
and they determine the geometry of the Seifert manifold (which could have different
fibrations) according to Table~\ref{tabellina}. The two non-Seifert
geometries are the Sol and the hyperbolic ones~\cite{Sco}.

\begin{table}
\begin{center}
\begin{tabular}{c|ccc} 
\phantom{\Big|} & $\chi^{\rm orb}>0$ & $\chi^{\rm orb}=0$ & $\chi^{\rm orb}<0$ \\ \hline
\phantom{\Big|} $e=0$ & $S^2\times\matR$ & $E^3$ & $H^2\times\matR$ \\
\phantom{\Big|} $e\neq 0$ & $S^3$ & Nil & $\widetilde{{\rm SL}_2\matR}$ \\ 
\end{tabular}
\end{center}
\caption{The six Seifert geometries.}
\label{tabellina}
\end{table} 

The following result shows how to compute the complexity (when $c\leqslant 10$)
of most manifolds belonging to the first 7 geometries. It is proved for $c\leqslant 9$
in~\cite{MaPe:geometric}, and completed for $c=10$ here in Subsection~\ref{bricks:found:subsection}.
We define the norm $|p,q|$ of two coprime non-negative integers inductively by setting
$|1,0|=|0,1|=|1,1|=0$ and $|p+q,q|=|p,q+p|=|p,q|+1$.
A norm $\|A\|$ on matrices $A\in\GL_2(\matZ)$ is
also defined in~\cite{MaPe:geometric}.

\begin{theorem}\label{non:hyperbolic:teo}
Let $M$ be a geometric non-hyperbolic manifold with $c(M)\leqslant 10$:
\begin{enumerate}
\item if $M$ is a lens space $L(p,q)$, then $c(M)=|p,q|-2$;
\item if $M$ is a torus bundle with monodromy $A$
then $c(M)=\min\{\|A\|+5, 6\}.$
\item \label{23m} if $M=\big(S^2,(2,1),(3,1),(m,1),-1\big)$ with $m\geqslant 5$, we have $c(M)=m$;
\item \label{2nm} if $M=\big(S^2,(2,1),(n,1),(m,1),-1\big)$ is not of the type above, 
  we have $c(M)=n+m-2$;
\item \label{23pq} if $M=\big(S^2,(2,1),(3,1),(p,q),-1\big)$ with $p/q>5$ is not of the types above, 
  we have $c(M)=|p,q|+2$;
\item if $M=\big(F, (p_1, q_1),\ldots, (p_k, q_k), t\big)$ is not of the types above, then
$$c(M)=\max\big\{0,t-1+\chi(F)\big\}+6\big(1-\chi(F)\big) +
  \sum_{i=1}^k\big(|p_i,q_i|+2\big). $$
\end{enumerate}
\end{theorem}

Note from Table~\ref{closed:table} that a Seifert manifold with $c<6$ has $\chi^{\rm orb}>0$ and
one with $c\leqslant 6$ has $\chi^{\rm orb}\geqslant 0$, 
whereas for higher $c$ most Seifert manifolds have $\chi^{\rm orb}<0$.

\begin{remark}
Theorem~\ref{non:hyperbolic:teo}, together with analogous formulas for
some non-geometric graph manifolds, follows from the decomposition of closed manifolds
into bricks, introduced in Section~\ref{bricks:section}.
The lists of all non-hyperbolic manifolds with $c\leqslant 10$ is then computed from
such formulas by a computer program, available from~\cite{weblist}.
A mistake in that program produced in~\cite{MaPe} for $c=9$ a list of $1156$ manifolds
instead of $1154$
(two graph manifolds with distinct parameters were counted twice). Using Turaev-Viro
invariants, Matveev has also recently checked that all the listed closed manifolds
with $c\leqslant 10$ are distinct~\cite{Mat:priv}.
\end{remark}

\subsection{Hyperbolic manifolds}
\begin{figure}[t]
\begin{minipage}{.03\textwidth}\hfil
\end{minipage}
\begin{minipage}{.2\textwidth}
\centering
\includegraphics[width = .9\textwidth]{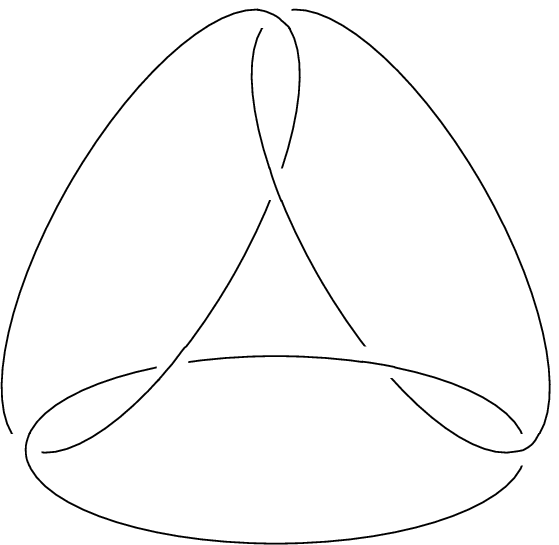}
\vglue -4mm
\end{minipage}
\begin{minipage}{.03\textwidth}\hfil
\end{minipage}
\begin{minipage}{.64\textwidth}
{The symmetries of this link 
act transitively on the components, in such a way 
that to define the $(p/q, r/s, t/u)$-surgery
we do not need to associate a component to each parameter.
}
\end{minipage}
\caption{The chain link with $3$ components.}
\label{chainlink:fig}
\end{figure}

Table~\ref{hyperbolic:table} shows all closed hyperbolic manifolds with $c\leqslant 10$.
Each such manifold is a Dehn surgery on the chain link with 3 components shown in
Fig.~\ref{chainlink:fig}, with parameters
shown in the table. 

\begin{table}
  \begin{center}
    \begin{tabular}{ccccc}
      surgery parameters & N. & volume & shortest geod & homology \\
      \multicolumn{4}{c}{complexity $9$\phantom{\Big|}} \\	
      $1,-4,-3/2$     & $1$  & $0.942707362$ & $0.5846$ & $\matZ_5 + \matZ_5$ \\ 
      $1,-4,2$        & $2$  & $0.981368828$ & $0.5780$ & $\matZ_5$ \\ 
      $1,-5,-1/2$     & $3$  & $1.014941606$ & $0.8314$ & $\matZ_3 + \matZ_6$ \\ 
      $1,-3/2,-3/2$   & $4$  & $1.263709238$ & $0.5750$ & $\matZ_5 + \matZ_5$ \\
      \multicolumn{4}{c}{complexity $10$\phantom{\Big|}} \\	
      $1,-5,2$        & $5$  & $1.284485300$ & $0.4803$ & $\matZ_6$ \\
      $1,2,1/2$       & $6$  & $1.398508884$ & $0.3661$ & trivial \\
      $1,-5,1/2$      & $7$  & $1.414061044$ & $0.7941$ & $\matZ_6$ \\
      $1,-4,3$        & $8$  & $1.414061044$ & $0.3648$ & $\matZ_{10}$ \\
      $1,-4,-4/3$     & $9$  & $1.423611900$ & $0.3523$ & $\matZ_{35}$ \\
      $1,2,-1/2$      & $10$ & $1.440699006$ & $0.3615$ & $\matZ_3$ \\
      $1,2,-3/2$      & $12$ & $1.529477329$ & $0.3359$ & $\matZ_5$ \\
      $1,-4,-5/2$     & $13$ & $1.543568911$ & $0.3353$ & $\matZ_{35}$ \\
      $1,-1/2,-5/2$   & $14$ & $1.543568911$ & $0.5780$ & $\matZ_{21}$ \\
      $1,-4,-5/3$     & $16$ & $1.583166660$ & $0.2788$ & $\matZ_{40}$ \\
      $1,-6,-1/2$     & $17$ & $1.583166660$ & $0.5577$ & $\matZ_{21}$ \\
      $1, -1/2, -7/2$ & $18$ & $1.583166660$ & $0.7774$ & $\matZ_3 + \matZ_9$ \\
      $2,-3/2,-3/2$   & $19$ & $1.588646639$ & $0.3046$ & $\matZ_{30}$ \\
      $1,-5,-3/2$     & $20$ & $1.588646639$ & $0.5345$ & $\matZ_{30}$ \\
      $1,-4,3/2$      & $21$ & $1.610469711$ & $0.2499$ & $\matZ_5$ \\
      $1,2,-5/2$      & $24$ & $1.649609715$ & $0.2627$ & $\matZ_7$ \\
      $1,-1/2,-3/2$   & $25$ & $1.649609715$ & $0.5087$ & $\matZ_{15}$ \\
      $1,1/2,-6$      & $34$ & $1.757126029$ & $0.7053$ & $\matZ_7$ \\
      $1,-1/2,-1/2$   & $49$ & $1.824344322$ & $0.4680$ & $\matZ_3 + \matZ_3$ \\
      $1,-5,-1/3$     & $55$ & $1.831931188$ & $0.5306$ & $\matZ_2 + \matZ_{12}$ \\
      $1,-3/2,-5/3$   & $74$ & $1.885414725$ & $0.3970$ & $\matZ_{40}$ \\
      $1,-5/2,-5/2$   & $76$ & $1.885414725$ & $0.5846$ & $\matZ_7 + \matZ_7$ \\
      $-5/2,-1/2,-1/2$& $77$ & $1.885414725$ & $0.5846$ & $\matZ_{39}$ \\
      $1,-5,-2/3$     & $91$ & $1.910843793$ & $0.4421$ & $\matZ_{30}$ \\
      $1,-4/3,-3/2$   & $139$ & $1.953708315$ & $0.3535$ & $\matZ_{35}$
		\end{tabular}
	\end{center}
  \caption{The hyperbolic manifolds of complexity $9$ and $10$. Each such manifold is
	described as the surgery on the chain link with some parameters.}
  \label{hyperbolic:table}
\end{table}

It is proved in~\cite{MaFo} that every closed 3-manifold with $c\leqslant 8$
is a graph manifold, and that the first closed hyperbolic manifolds arise with $c=9$. The hyperbolic
manifolds with $c=9$ then turned out~\cite{MaPe} to be the 4 smallest ones known. 
The most interesting question about those with $c=10$ is then whether they are also
among the smallest ones known, for instance
comparing them with the closed census~\cite{HoWe} also used by SnapPea~\cite{SnapPea}.
As explained in~\cite{DuTh}, the manifolds in that census have all geodesics bigger than $.3$,
and therefore some manifolds having $c=10$ are not present there
(namely, those in Table~\ref{hyperbolic:table} corresponding to N.~$16$, $21$, $24$).
We have therefore used SnapPea (in the python version)
to compute a list of many surgeries on the chain link with $3$ components 
(avoiding the non-hyperbolic ones, listed in~\cite{MaPe:chain}),
available from~\cite{weblist}, which contains many closed manifolds of volume smaller than $2$
that are not present in SnapPea's closed census. The entry ``N.'' in Table~\ref{hyperbolic:table}
tells the position of the manifold in our table from~\cite{weblist}.
The first $10$ manifolds of the two
lists nevertheless coincide and are also fully described in~\cite{HoWe}, and
they all have $c\leqslant 10$, as Table~\ref{hyperbolic:table} shows.

\subsection{Non-geometric manifolds}
Every non-hyperbolic orientable manifold with $c\leqslant 10$ is a graph manifold,
\emph{i.e.}~its \JSJ\ decomposition consists of Seifert or Sol blocks. 
A non-geometric orientable 
manifold whose decomposition contains a hyperbolic block with $c\leqslant 11$
is constructed in~\cite{AmMa}, and from our census now it follows that it cannot
have $c\leqslant 10$. Therefore we have proved the following.
\begin{theorem} The first closed orientable irreducible manifold with non-trivial \JSJ\ decomposition
containing hyperbolic blocks has $c=11$.
\end{theorem}
All graph manifolds with $c\leqslant 10$ are collected in Table~\ref{Seifert:table} according to 
their \JSJ\ decomposition into fibering pieces, and to the type of fiberings of each piece.
\begin{table}
  \begin{center}
    \begin{tabular}{rccccccccccc}
      & $0$ & $1$ & $2$ & $3$ & $4$ & $5$ & $6$ & $7$ & $8$ & $9$ & $10$ \\
      
      \multicolumn{12}{c}{geometric\phantom{\Big|}} \\

      lens spaces &
      $3$ & $2$ & $3$ & $6$ & $10$ & $20$ & $36$ & $72$ & $136$ & $272$ & $528$ \\

			$S^2, 3$
			& . & . & $1$ & $1$ & $4$ & $11$ & $31^*$ & $84$ & $226$ & $586$ & $1477$ \\

			$S^2, 4$ 
			& . & . & .   & .   & .   & .    & $2$    & $4$  & $14$  & $40$  & $120$ \\

			$S^2, 5$ 
			& . & . & .   & .   & .   & .    & .      & .    & .     & $2$   & $5$ \\

			$\matRP^2, 2$ 
			& . & . & .   & .   & .   & .    & $2$    & $4$  & $14$  & $34$  & $90$ \\

			$\matRP^2, 3$ 
			& . & . & .   & .   & .   & .    & .      & .    & .     & $2$   & $5$ \\

			$T$ or $K$  
			& . & . & .   & .   & .   & .    & $4^*$  & $2$  & $2$   & $2$   & $2$ \\

			$T,1$ or $K,1$  
			& . & . & .   & .   & .   & .    & .      & .    & .     & $4$   & $10$\\

			$T$-fiberings over $S^1$  
			& . & . & .   & .   & .   & .    & .      & $2$  & $2$   & $6$   & $6$\\

			$T$-fiberings over $I$ 
			& . & . & .   & .   & .   & .    & .      & $3$  & $7$   & $17$  & $33$\\

      \multicolumn{12}{c}{non-geometric\phantom{\Big|}} \\

			$D, 2$ --- $D, 2$
			& . & . & .   & .   & .   & .    & .      & $4$  & $35$  & $168$ & $674$\\

			$A, 1$ 
			& . & . & .   & .   & .   & .    & .      & .    & .     & $8$   & $24$ \\

			$D, 2$ --- $D, 3$ 
			& . & . & .   & .   & .   & .    & .      & .    & .     & $3$   & $24$ \\

			$S, 1$ --- $D, 2$  
			& . & . & .   & .   & .   & .    & .      & .    & .     & $3$   & $24$ \\

			$D, 2$ --- $A, 1$ --- $D, 2$ 
			& . & . & .   & .   & .   & .    & .      & .    & .     & $3$   & $31$ \\

      total  
			& $\bf{3}$ & $\bf{2}$ & $\bf{4}$ & $\bf{7}$ & $\bf{14}$ & $\bf{31}$ & $\bf{74}$ & $\bf{175}$ & $\bf{436}$ & $\bf{1150}$ & $\bf{3053}$ \\

    \end{tabular}
  \end{center} 
  \caption{The type of graph manifolds of
    given complexity, up to $10$. Here, $I, D, S, A, T, K$ denote respectively the closed interval,
		the disc, the M\"obius strip, the annulus, the torus, and the Klein bottle. We denote by
		$X,n$ a block with base space the surface $X$ and $n$ exceptional fibers. We write $X$ for
		$X,0$. We have counted as $T$-fiberings only the Sol manifolds, not the manifolds also admitting
		a Seifert structure. 
		There is a flat manifold with $c=6$ counted twice, since it has two different fibrations,
		corresponding to the asterisks.}
  \label{Seifert:table}
\end{table}

\subsection{The simplest manifolds}
As the following discussion shows, in most geometries,
the manifolds with lowest complexity are the ``simplest'' ones.

\subsubsection{Elliptic} The elliptic manifolds of smallest complexity are 
$S^3, \matRP^3,$ and $L(3,1)$, having $c=0$.
The first manifold which is not a lens space is
$\big(S^2,(2,1),(2,1),(2,1),-1\big)$ and has $c=2$. 
It is the elliptic manifold with
smallest non-cyclic fundamental group, having order $8$~\cite{Mat:book}. 

\subsubsection{Flat}
Every (orientable or not) flat manifold has $c=6$. A typical way to obtain some flat
3-manifold $M$ is from a face-pairing of the cube: by taking a triangulation
of the cube with $6$ tetrahedra matching along the face-pairing, 
we get a minimal triangulation of $M$.

\subsubsection{$\matH^2\times\matR$}
The first manifolds of type $\matH^2\times\matR$ are non-orientable and have $c=7$, and are
also the manifolds of that geometry with smallest base orbifold~\cite{AmMa2}, 
having volume 
$-2\pi\chiorb = \pi/3$. 

\subsubsection{Sol}
The first manifold of type Sol 
is also non-orientable and has $c=6$,
and it is the unique filling of the Gieseking manifold, the cusped 
hyperbolic manifold with smallest
volume $1.0149\ldots$~\cite{Ad} and smallest 
complexity $1$~\cite{CaHiWe}.
It is also the unique torus
fibering whose monodromy $A=\matr 0111$ is hyperbolic with $|\tr A|<2$~\cite{AmMa2}.

\subsubsection{Hyperbolic}
As we said above, the first orientable 
hyperbolic manifolds are the smallest ones known. It would be
interesting to know the complexity of the first non-orientable closed hyperbolic manifold,
whose volume is probably considerably bigger than in the orientable case, see~\cite{HoWe}.

\section{Census of hyperbolic manifolds} \label{hyperbolic:section}
We describe here the compact hyperbolic manifolds with boundary with $\chi=0$ and $c\leqslant 7$,
and the orientable ones with $\chi<0$ and $c\leqslant 4$.

\subsection{Manifolds with $\chi = 0$}
Recall that we define a compact $M$ to be
hyperbolic when it admits a complete metric
of finite volume and geodesic boundary, after removing 
all boundary components with $\chi = 0$.
Therefore, hyperbolic manifolds $M$ with $\chi(M)=0$ have some cusps based on tori or Klein
bottles, and those with $\chi(M)<0$ have geodesic boundary and possibly some cusps.
To avoid confusion, we define the \emph{topological boundary} of $M$ 
to be the union of the geodesic boundary and the cusps.

Hyperbolic manifolds with $\chi(M)=0$ and $c\leqslant 7$ were listed 
by Hodgson and Weeks in~\cite{CaHiWe}
and form the cusped census used by SnapPea. They are collected, according to their
topological boundary, in Table~\ref{cusped:table}.
Hyperbolicity of each manifold was checked by
solving Thurston's equations, and all manifolds were distinguished computing their
Epstein-Penner \emph{canonical decomposition}~\cite{EpPe}. In practice,
volume, homology, and 
the length of the shortest geodesic are usually enough to distinguish two such manifolds.

\begin{table}
  \begin{center}
    \begin{tabular}{rcccccccc}
      topological boundary & $0$ & $1$ & $2$ & $3$ & $4$ & $5$ & $6$ & $7$ \\
      
      \multicolumn{9}{c}{orientable\phantom{\Big|}} \\

      $T$ & . & . & $2$ & $9$ & $52$ & $223$ & $913$ & $3388$ \\
      
      $T,T$ & . & . & .   & .   & $4$  & $11$  & $48$  & $162$ \\
      
      $T,T,T$ & . & . & .   & .   & .    & .     & $1$   & $2$ \\

      total orientable & . & . & $\bf 2$ & $\bf 9$ & $\bf{56}$ & $\bf{234}$ & $\bf{962}$ & $\bf{3552}$ \\
      
      \multicolumn{9}{c}{non-orientable\phantom{\Big|}} \\

      $K$ & . & $1$ & $1$ & $5$ & $14$ & $52$ & $171$ & $617$ \\
      
      $K,K$ & . & .   & $1$ & $2$ & $9$  & $23$ & $68$  & $208$ \\
      
      $K,K,K$ & . & .   & .   & .   & .    & .    & $3$   & $6$   \\
      
      $K,K,K,K$ & . & .   & .   & .   & .    & .    & $1$   & .     \\
      
      $T$ & . & .   & .   & .   & $1$  & $1$  & $4$   & $19$  \\
      
      $T,T$ & . & .   & .   & .   & .    & .    & $1$   & .     \\
      
      $K,T$ & . & .   & .   & .   & $1$  & $2$  & $8$   & $31$  \\
      
      $K,K,T$ & . & .   & .   & .   & $1$  & .    & $3$   & $6$   \\

      total non-orientable 
      & . & $\bf 1$ & $\bf 2$ & $\bf 7$ & $\bf{26}$ & $\bf{78}$ & $\bf{259}$ & $\bf{887}$

    \end{tabular}
  \end{center}
  \caption{The number of cusped hyperbolic manifolds of
    given complexity, up to $7$. The ``topological boundary'' indicates the tori $T$ and Klein bottles
    $K$ present as cusps.}
  \label{cusped:table}
\end{table}

\subsection{Manifolds with $\chi<0$} \label{Kojima:subsection}
Equations analogous to Thurston's were constructed by Frigerio and Petronio in~\cite{FriPe}
for an ideal triangulation $T$ of a manifold $M$ with $\chi(M)<0$.
A solution of such equations gives a realization of the hyperbolic structure of $M$
via partially truncated hyperbolic tetrahedra. One such tetrahedron is
parametrized by its $6$ interior dihedral angles $\alpha_1,\ldots,\alpha_6$.
The sum of the $3$ of them incident to a given vertex must be less or equal than $\pi$,
and the vertex is truncated if the sum is less than $\pi$, or ideal if it is $\pi$.
The compatibility equations ensure that identified edges all have the same length and
that dihedral angles sum to $2\pi$ around each resulting edge. These equations, together with others 
checking the completeness
of the cusps, realize the hyperbolic structure for $M$. Then Kojima's \emph{canonical
decomposition}~\cite{Koj}, analogous to Epstein-Penner's,
is a complete invariant which allows one 
to distinguish manifolds. In contrast with the case $\chi=0$, there are plenty of manifolds
having the same complexity that are not distinguished by volume, homology, Turaev-Viro invariants,
and the canonical decomposition seems to be the only available tool, see 
Subsection~\ref{families:subsection}.
The results from~\cite{FriMaPe2} are summarized in Table~\ref{geodesic:table}.

\begin{table}
  \begin{center}
    \begin{tabular}{rccccc}
      topological boundary 
      & $0$ & $1$ & $2$ & $3$ & $4$ \\
      
      $2$ & . & . & $8$ & $76$ & $628$ \\
      
      $3$ & . & . & .    & $74$ & $2034$ \\
      
      $4$ & . & . & .    & .    & $2340$ \\
      
      $2,0$ & . & . & .    & $1$  & $18$ \\
      
      $3,0$ & . & . & .    & .    & $12$ \\
      
      $2,0,0$ & . & . & .    & .    & $1$ \\ 
      
      total   & . & . & $\bf 8$ & $\bf{151}$ & $\bf{5033}$ \\
      
    \end{tabular}
  \end{center}
  \caption{The number of orientable hyperbolic manifolds with non-empty geodesic boundary of
    given complexity, up to $4$. The ``topological boundary'' indicates the genera of the boundary
    components, with zeroes correspond to cusps.}
  \label{geodesic:table}
\end{table}

\begin{remark}
The two censuses of hyperbolic manifolds described in this Section have a slightly more
experimental nature than the closed census of Section~\ref{closed:section},
since solving hyperbolicity equations and calculating the canonical decomposition
involve numerical calculations with truncated digits.
\end{remark}

\section{Complexity and volume of hyperbolic manifolds} \label{complexity:volume:section}
We describe here some relations between the complexity and the volume
of a hyperbolic 3-manifold.

\subsection{Ideal tetrahedra and octahedra}
As Theorem~\ref{vol:c:teo} below shows, there is a constant $K$ such that
${\rm Vol} (M)< K \cdot c(M)$ for any hyperbolic $M$.
Let $v_{\rm T} = 1.0149\ldots $ and $v_{\rm O} = 3.6638\ldots$ be the volumes 
respectively of the regular ideal hyperbolic tetrahedron and octahedron.
\begin{theorem} \label{vol:c:teo}
Let $M$ be hyperbolic, with or without boundary. 
If $\chi(M)=0$ we have ${\rm Vol}(M) \leqslant v_{\rm T}\cdot c(M)$. 
If $\chi(M)<0$ we have ${\rm Vol} (M)< v_{\rm O}\cdot c(M)$.
\end{theorem}
\begin{proof}
First, note that by the naturality property of the complexity $c(M)$ is
the minimum number of tetrahedra in an (ideal) triangulation.
If $M$ is closed, take a minimal triangulation $T$ and straighten it. 
Tetrahedra may overlap or collapse to low-dimensional objects, having
volume zero. Since geodesic tetrahedra have volume less than $\voltet$,
we get the inequality.

If $M$ is not closed, let $T$ be an ideal triangulation for $M$ 
with $c(M)$ tetrahedra. 
We can realize topologically $M$ with its boundary tori removed, 
by partially truncating
each tetrahedron in $T$ (\emph{i.e.}~removing the vertex only in presence of 
a cusp, and an open
star of it in presence of true boundary). Then we can
straighten every truncated tetrahedron with respect to the hyperbolic structure in $M$.
As above, tetrahedra may overlap or collapse. 
In any case,
the volume of each such will be at most $v_{\rm T}$ if there is no boundary, and strictly less than
$v_{\rm O}$ in general, since any ideal tetrahedron has volume at most equal to $v_{\rm T}$,
and any partially truncated tetrahedron has volume strictly less than $v_{\rm O}$~\cite{Ush}.
\end{proof}

The constants $v_{\rm T}$ and $v_{\rm O}$ are the best possible ones, see 
Remark~\ref{best:constants:rem}.
A converse result of type $c(M) < K' \cdot {\rm Vol} (M)$ is impossible,
because for big $C$'s there are a finite number of hyperbolic manifolds 
with complexity
less than $C$, and an infinite number of such with volume less than $C$.

\subsection{First segments of $c$ and $\Vol$}
Complexity and volume give two partial orderings on the set $\calH$ of all hyperbolic
3-manifolds. By what was just said, they are globally qualitatively very different.
Nevertheless, as noted in~\cite{MaFo}, they might have similar behaviours
on some subsets of $\calH$. We propose the following conjecture.

\begin{conjecture} \label{volume:complexity:conj}
Among hyperbolic manifolds with the same topological boundary, the ones with
smallest complexity have volume smaller than the other ones.
\end{conjecture}
The conjecture is stated more precisely as follows:
let $\calM_{\Sigma}$ be the set of hyperbolic manifolds having some fixed topological
boundary $\Sigma$. Suppose $M\in\calM_{\Sigma}$ is so that $c(M')\geqslant c(M)$ for all
$M'\in\calM_{\Sigma}$. We conjecture that $\Vol(M')>\Vol(M)$
for all $M'\in\calM_{\Sigma}$ having $c(M')>c(M)$.
We now discuss our conjecture.

\subsubsection{Closed case}
The closed hyperbolic manifolds with smallest $c=9$ are the four having smallest
volume known, see Table~\ref{hyperbolic:table}. Therefore
Conjecture~\ref{volume:complexity:conj} claims that these four are actually the ones
having smallest volumes among all closed hyperbolic manifolds.

\subsubsection{Connected topological boundary}
In this case, Conjecture~\ref{volume:complexity:conj} is true, as the following
shows.
\begin{theorem} \label{connected:surface:teo}
Among hyperbolic manifolds whose topological boundary is a connected surface,
the ones with smallest volume are the ones with smallest complexity.
\end{theorem}
\begin{proof}
Among manifolds having one toric cusp,
the figure-8 knot complement and its sibling are those with minimal volume
$2\cdot\voltet$~\cite{CaMe} and minimal complexity $2$. Among those with a Klein bottle cusp,
the Gieseking manifold is the one with minimal volume $\voltet$~\cite{Ad} 
and minimal complexity $1$. Our assertion restricted to orientable 3-manifolds
bounded by a connected surface of higher genus is proved in~\cite{FriMaPe} combining
the naturality property of the complexity with
Miyamoto's description~\cite{Miy} of all such manifolds with minimal volume. 
The same proof also works in the general case.
\end{proof}

\subsubsection{Experimental data}
Conjecture~\ref{volume:complexity:conj} is true
when restricted to the manifolds of Tables~\ref{hyperbolic:table},
~\ref{cusped:table}, and~\ref{geodesic:table}, for all the boundary types involved
(see~\cite{CaHiWe},~\cite{SnapPea}, and~\cite{FriMaPe2}).
One sees from Table~\ref{hyperbolic:table} that the manifolds of type
$(K,K)$, $(T,T)$, $(K,T)$, $(K,K,T)$, $(T,T,T)$, $(K,K,K)$, and $(K,K,K,K)$ 
with smallest complexity have respectively $c=2,4,4,4,6,6,$ 
and $6$. The manifolds with $c=2$ are constructed with two regular ideal tetrahedra,
and hence have volume $2\cdot\voltet$. Those with $c=4$ are constructed either with
$4$ regular ideal tetrahedra, hence having volume $4\cdot\voltet = 4.05976\ldots$,
or with one regular ideal octahedron, of volume $\voloct = 3.6638\ldots$ (therefore
Conjecture~\ref{volume:complexity:conj} claims that every other $M$ with the same
topological boundary has volume bigger than $4\cdot\voltet$). Those with $c=6$
have volume $2\cdot\voldrum = 5.3334\ldots$, where  $\voldrum = 2.6667\ldots$ is the volume of 
the ``triangular ideal drum''
used by Thurston~\cite{Thu} to construct
the complement of the chain link of Fig.~\ref{chainlink:fig}, which is the only
orientable manifold among them.

\begin{problem}
Classify the hyperbolic (orientable) manifolds of smallest complexity among those having
$\chi=0$ and $k$ toric cusps, and compute their volume, for each $k$.
\end{problem}

\subsection{Matveev-Fomenko conjecture} \label{counterexample:subsection}
As we mentioned above, the orderings given by $c$ and $\Vol$ are qualitatively different
on the whole set $\calM$ of hyperbolic manifolds, but might be similar on some subsets
of $\calM$. The following conjecture was proposed by Matveev and Fomenko in~\cite{MaFo}.
\begin{conjecture}[Matveev-Fomenko~\cite{MaFo}] \label{MaFo:conj}
Let $M$ be a hyperbolic manifold with one cusp. 
Among Dehn fillings $N$ and $N'$ of $M$, if $c(N)<c(N')$ then
$\Vol(N)<\Vol(N')$.
\end{conjecture}
The complexity-$10$ closed census produces a counterexample to Conjecture~\ref{MaFo:conj}.
\begin{proposition}
Let $N(p/q)$ be the $p/q$-surgery on the figure-8 knot. We have
\begin{eqnarray*}
\Vol\big(N(5/2)\big) = 1.5294773\ldots & \quad & c\big(N(5/2)\big) = 11 \\
\Vol\big(N(7)\big) =   1.4637766\ldots & \quad & c\big(N(7)\big) > 11
\end{eqnarray*}
\end{proposition}
\begin{proof}
We first note that $N(p/q) = N(-p/q)$ is the $(1,2,1-p/q)$-surgery on the chain link.
The manifold $N(7)$ does not belong to Table~\ref{hyperbolic:table} (it is the manifold
labeled as N.11 in our census of surgeries on the chain link of~\cite{weblist}), and
hence has $c>11$, whereas $N(5/2)$ is the manifold N.12 and has $c=11$.
\end{proof}

\section{Lower bounds} \label{lower:bounds:section}
Providing upper bounds for the complexity of a given manifold $M$ is relatively easy:
from any combinatorial description of $M$ one recovers a spine of $M$ with $n$ 
vertices, and certainly $c(M)\leqslant n$. 
Finding lower bounds is a much more difficult task. 
The only $\partial$-irreducible
manifolds whose complexity is known are those listed in the censuses
of Sections~\ref{closed:section} and~\ref{hyperbolic:section}, and some infinite families
of hyperbolic manifolds with bundary described below.
In particular, for a closed irreducible $M$, the value
$c(M)$ is only known when $c(M)\leqslant 10$, \emph{i.e.} for a finite number of manifolds.
\subsection{The closed case}
The only available lower bound 
for closed irreducible orientable manifolds is the following one, due to Matveev and 
Pervova. We denote by $|{\rm Tor}(H_1(M))|$ the order of the torsion subgroup
of $H_1(M)$, while $b_1$ is the rank of the free part, 
\emph{i.e.}~the fist Betti number of $M$.
\begin{theorem}[Matveev-Pervova~\cite{MatPer}]
Let $M$ be a closed orientable irreducible manifold different from $L(3,1)$. 
Then $c(M)\geqslant
2\cdot\log_5|{\rm Tor}(H_1(M))|+b_1-1$.
\end{theorem}
Recall that Theorem~\ref{non:hyperbolic:teo} holds only for $c\leqslant 10$. Actually,
the same formulas in the statement give an upper bound for $c(M)$. Some such
upper bounds for lens spaces, torus bundles, and simple Seifert manifolds
were previously found by Matveev and Anisov, 
who proposed the following conjectures.
\begin{conjecture}[Matveev~\cite{Mat:book}] We have
$$c\big(L(p,q)\big) = |p,q|-2 \quad {\rm and} \quad c\big(S^2,(2,1),(2,1),(m,1),-1\big) = m$$
\end{conjecture}
\begin{conjecture}[Anisov~\cite{An:towards}]
The complexity of a torus bundle $M$ over $S^1$ with monodromy $A\in\GL_2(\matZ)$ is
$c(M) = \min\{\|A\|+5, 6\}.$
\end{conjecture}

\subsection{Families of hyperbolic manifolds with boundary of known complexity} 
\label{families:subsection}
The following corollaries of Theorem~\ref{vol:c:teo} were first noted by Anisov.
\begin{corollary}[Anisov~\cite{An}]
The complexity of a hyperbolic manifold decomposing into $n$ ideal regular tetrahedra
is $n$.
\end{corollary}
\begin{corollary}[Anisov~\cite{An}]
The punctured torus bundle with monodromy $\matr 2111^n$ is a hyperbolic manifold
of complexity $2n$.
\end{corollary}
For each $n\geqslant 2$, 
Frigerio, Martelli, and Petronio defined~\cite{FriMaPe} the family $\calM_n$  
of all orientable compact manifolds admitting an ideal triangulation with
one edge and $n$ tetrahedra.
\begin{theorem}[Frigerio-Martelli-Petronio~\cite{FriMaPe}] \label{FriMaPe:teo}
Let $M\in\calM_n$. Then $M$ is hyperbolic with a genus-$n$ surface as geodesic boundary, and
without cusps. It has complexity $n$. Its homology, volume, 
Heegaard genus, and Turaev-Viro invariants also depend only on $n$. 
\end{theorem}
The manifolds in $\calM_n$ are distinguished by their Kojima's canonical decomposition 
(see Subsection~\ref{Kojima:subsection}), which is precisely the triangulation with one edge
defining them. Therefore combinatorially different such triangulations give different
manifolds. 
\begin{theorem}[Frigerio-Martelli-Petronio~\cite{FriMaPe, FriMaPe3}] \label{FriMaPe:growth:teo}
Manifolds in $\calM_n$ correspond bijectively to triangulations with one edge and $n$ tetrahedra.
The cardinality $\#\calM_n$ grows as $n^n$.
\end{theorem}
We say that a sequence $a_n$ \emph{grows as $n^n$} when 
there exist constants $0<k<K$ such that $n^{k\cdot n}<a_n<n^{K\cdot n}$ for all $n\gg 0$.
\begin{corollary}[Frigerio-Martelli-Petronio~\cite{FriMaPe3}] \label{hyp:growth:teo}
The number of hyperbolic manifolds of complexity $n$ grows as $n^n$.
\end{corollary}
\begin{remark} \label{best:constants:rem}
From the families introduced here we see that the inequalities of Theorem~\ref{vol:c:teo}
cannot be strengthened. The torus bundles $M$ above have $\Vol(M)=\voltet\cdot c(M)$,
and the manifolds in $\calM_n$ have $\Vol(M) = v_n\cdot c(M)$, with $v_n$ equals to
the volume of a truncated tetrahedron with all angles $\pi/(3n)$, so that
$v_n\to\voloct$ for $n\to\infty$.
\end{remark}

The set $\calM_n$ is also the set mentioned in Theorem~\ref{connected:surface:teo}
of all manifolds having both minimal complexity and minimal volume among those with
a genus-$n$ surface as boundary. We therefore get from Table~\ref{geodesic:table} that
$\#\calM_n$ is $8, 74, 2340$ for $n=2,3,4$.

The class $\calM_n$ is actually contained as $\calM_n = \calM_{n,0}$ in a bigger family
$\calM_{g,k}$, defined in~\cite{FriMaPe3}. The set $\calM_{g,k}$ consists of 
all orientable hyperbolic manifolds of 
complexity $g+k$ with connected geodesic boundary of genus $g$ and $k$ cusps. 
Theorems~\ref{FriMaPe:teo} and~\ref{FriMaPe:growth:teo} hold similarly for all such sets.
For any fixed $g$ and $k$, $\calM_{g,k}$ is the set of all manifolds with minimum complexity
among those with that topological boundary. Therefore Conjecture~\ref{volume:complexity:conj}
would imply the following.
\begin{conjecture}[Frigerio-Martelli-Petronio~\cite{FriMaPe3}] The manifolds of smallest volume
among those with a genus-$g$ geodesic surface as boundary and $k$ cusps are those in
$\calM_{g,k}$.
\end{conjecture}

\section{Minimal spines} \label{minimal:section}
We describe here some known results about minimal spines, which are crucial for computing 
the censuses of Sections~\ref{closed:section} and~\ref{hyperbolic:section}.
\subsection{The algorithm} \label{algorithm:subsection}
The algorithm used to classify all manifolds with increasing complexity $n$ typically works
as follows: 
\begin{enumerate}
\item list all special spines with $n$ vertices (or triangulations with $n$ tetrahedra);
\item remove from the list the many spines that are easily seen to be 
non-minimal, or not to thicken to an irreducible
(or hyperbolic) manifold; 
\item try to recognize the manifolds obtained from thickening
the remaining spines;
\item eliminate from that list of manifolds the duplicates, and the manifolds
that have already been found previously in some complexity-$n'$ census for some $n'<n$.
\end{enumerate}
Typically, step (1) produces a huge list of spines, $99.99\ldots$ \% of which are canceled via
some quick criterion of non-minimality during step (2), and one is left with a much smaller
list, so that steps (3) and (4) can be done by hand.

\subsection{Cutting dead branches} \label{branches:subsection}
Step (1) of the algorithm above needs a huge
amount of computer time already for $c=5$, due to the
very big number of spines listed.
Therefore one actually uses the non-minimality criteria (step (2)) \emph{while} 
listing the special spines
with $n$ vertices (step (1)), to cut many ``dead branches''. Step (1) 
remains the most expensive one
in terms of computer time, so it is worth describing it with some details. 

A special spine or its dual (possibly ideal) 
triangulation $T$ (see Remark~\ref{spine:triangulation:rem}) with $n$ tetrahedra
can be encoded roughly as follows.
Take the face-pairing 4-valent graph $G$ of the tetrahedra in $T$. 
It has $n$ vertices and $2n$ edges.
After fixing a simplex on each vertex, a
label in $S_3$ on each (oriented) edge of $G$ encodes how the faces are glued.
We therefore get $6^{2n}$ gluings 
(the same combinatorial $T$ is usually realized by many distinct gluings).
Point (1) in the algorithm consists of two steps:
\begin{itemize}
\item[(1a)] classify all 4-valent graphs $G$ with $n$ vertices;
\item[(1b)] for each graph $G$, fix a simplex on each vertex, and 
try the $6^{2n}$ possible labelings on edges.
\end{itemize}
Step (1b) is by far the most expensive one, because it contains
many ``dead branches''; most of them are cut as follows: 
a partial labeling of some $k$ of the $2n$ edges defines a partial gluing of the tetrahedra.
If such partial gluing already fulfills some local non-minimality criterion,
we can forget about every labeling containing this partial one.

\begin{remark} \label{ografi:rem}
A spine of an \emph{orientable} manifold can be encoded more efficiently by fixing
an immersion of the graph $G$
in $\matR^2$, and assigning a colour in $\matZ_2$ to each vertex and a colour in
$\matZ_3$ to each edge~\cite{BePe}.
\end{remark}

Local non-minimality criteria used to cut the branches are listed in Subsection~\ref{criteria:subsection}.
We discuss in Subsection~\ref{graph:subsection} another powerful tool, which works
in the closed case only: it turns out that
most $4$-valent graphs $G$ can be quickly checked \emph{a priori} not to give rise to any
minimal spine (of closed manifolds). 
\subsection{Local non-minimality criteria} \label{criteria:subsection}
We start with the following results.
\begin{proposition}[Matveev~\cite{Mat90}] \label{small:faces:prop}
Let $P$ be a minimal special spine of a 3-manifold $M$. Then $P$ contains no embedded
face with at most $3$ edges.
\end{proposition}
\begin{proposition}[Matveev~\cite{Mat90}] \label{counterpass:prop}
Let $P$ be a minimal special spine of a closed orientable 3-manifold $M$.
Let $e$ be an edge of $P$.
A face $f$ cannot be incident $3$ times to $e$, and it cannot run twice on $e$ with opposite
directions.
\end{proposition}
In the orientable setting, both Propositions~\ref{small:faces:prop} and~\ref{counterpass:prop}
are special cases of the following. Recall that $S(P)$ is the subset of a special spine $P$
consisting of all points of type (2) and (3) shown in Fig.~\ref{special:fig}.
\begin{proposition}[Martelli-Petronio~\cite{MaPe}]
Let $P$ be a minimal spine of a closed orientable 3-manifold $M$. 
Every simple closed curve
$\gamma\subset P$ bounding a disc in the ball $M\setminus P$ and intersecting $S(P)$
transversely in at most $3$ points is contained in a small neighborhood of a point of $P$.
\end{proposition}
Analogous results in the possibly non-orientable setting are proved by Burton~\cite{Bu:graphs}.

\subsection{Four-valent graphs} \label{graph:subsection}
\begin{figure}
\centering
\includegraphics[width = 11 cm] {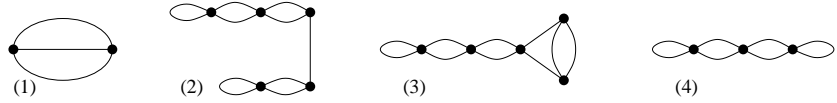}
\caption{Portions of graphs forbidden for minimal triangulations/spines of closed manifolds.}
\label{portions:fig}
\end{figure}

Quite surprisingly, some 4-valents graphs can be checked \emph{a priori} not to give
any minimal special spine of closed 3-manifold. 
\begin{remark}
The face-pairing graph of a (possibly ideal) triangulation is also
the set $S(P)$ in the dual special spine $P$.
\end{remark}
\begin{proposition}[Burton~\cite{Bu:graphs}] \label{Burton:prop}
The face-pairing graph $G$ of a minimal triangulation with at least $3$ tetrahedra
does not contain any portions
of the types shown in Fig.~\ref{portions:fig}-(1,2,3), except if $G$ itself is as in
Fig.~\ref{portions:fig}-(4). 
\end{proposition}
A portion of $G$ is of type shown in Fig.~\ref{portions:fig}-(2,3,4) when it 
is as in that picture, with chains of arbitrary length. 
In the algorithm of Subsection~\ref{branches:subsection}, 
step (1b) can be therefore restricted to the \emph{useful} 
4-valent graphs, \emph{i.e.}~the ones that do not
contain the portions forbidden by Proposition~\ref{Burton:prop}. Table~\ref{Burton:table},
taken from~\cite{Bu:graphs},
shows that some 40 \% of the graphs are useful.
\begin{table}
\begin{center}
\begin{tabular}{cccccccccc} 
\phantom{\Big|} $n$    & $3$ & $4$  & $5$  & $6$  & $7$   & $8$    & $9$    & $10$    & $11$ \\
\hline
\phantom{\Big|} useful & $2$ & $4$  & $12$ & $39$ & $138$ & $638$  & $3366$ & $20751$ & $143829$ \\
                all    & $4$ & $10$ & $28$ & $97$ & $359$ & $1635$ & $8296$ & $48432$ & $316520$
\end{tabular}
\end{center}
\caption{Useful graphs among all 4-valent graphs with $n\leqslant 11$ vertices.}
\label{Burton:table}
\end{table} 

\section{Bricks} \label{bricks:section}
As shown in Sections~\ref{definitions:section} and~\ref{minimal:section}, 
classifying all closed \ptwoirred\ 
manifolds with complexity $n$ reduces to listing all minimal
special spines of such manifolds with $n$ vertices. Non-minimality criteria as
those listed in Section~\ref{minimal:section} are then crucial to eliminate
the many non-minimal spines (by cutting ``dead branches'') 
and gain a lot of computer time. Actually, closed manifolds often have many minimal spines, and it is
not necessary to list them all: a criterion that eliminates some, but not all, minimal
spines of the same manifold is also suitable for us. This is the basic idea which underlies
the decomposition of closed \ptwoirred\ manifolds into \emph{bricks}, introduced by
Martelli and Petronio in~\cite{MaPe}, and described in the orientable case in this Section.
(For the nonorientable one, see~\cite{MaPe:nonori}.)

\subsection{A quick introduction}
The theory is roughly described as follows: every closed irreducible manifold $M$
decomposes along tori into pieces on which the complexity is additive. 
Each torus is marked with a $\theta$-graph in it, and the complexity of each piece is not
the usual one, because it depends on that graphs. A manifold $M$ which does not decompose is
a \emph{brick}. Every closed irreducible manifold decomposes into bricks. 
The decomposition is not unique, but there can be only a finite number of such.
In order to classify all manifolds with $c\leqslant 10$, one classifies all
bricks with $c\leqslant 10$, and then assemble them in all possible (finite) ways
to recover the manifolds. 

For $c\leqslant 10$, bricks are atoroidal, hence either Seifert or hyperbolic.
And the decomposition into bricks is tipically 
a mixure of the \JSJ, the graph-manifolds decomposition, and 
the thick-thin decomposition for hyperbolic manifolds. 
Very few closed manifolds do not decompose,
\emph{i.e.}~are themselves bricks.

\begin{proposition}
There are $25$ closed bricks with $c\leqslant 10$. They are: $24$ Seifert manifolds
of type $\big(S^2,(2,1),(m,1),(n,1),-1\big)$, and the hyperbolic manifold N.34 of
Table~\ref{hyperbolic:table}.
\end{proposition}
Among closed bricks, we have 
Poincar\'e homology sphere $\big(S^2,(2,1),(3,1),(5,1),-1\big).$
\begin{proposition}
There are $21$ non-closed bricks with $c\leqslant 10$. 
\end{proposition}
There are $4978$ closed irreducible manifolds with $c\leqslant 10$, 
see Table~\ref{closed:table}. Therefore $4953 = 4978-25$ such manifolds
are obtained with the $21$ bricks above.

Before giving precise definitions, we note that 
the \emph{layered triangulations}~\cite{Bu, JaRu} 
of the solid torus $H$ are particular decompositions of $H$
into bricks. Our experimental results show the following.
\begin{proposition}
Every closed irreducible atoroidal manifold with $c\leqslant 10$ has a minimal
triangulation containing a (possibly degenerate~\cite{Bu})
layered triangulation, except for some 
$\big(S^2,(2,1),(m,1),(n,1),-1\big)$ 
and the hyperbolic N.34 of Table~\ref{hyperbolic:table}.
\end{proposition}

\subsection{$\theta$-graphs in the torus}
In this paper, a \emph{$\theta$-graph} $\theta$ in the torus $T$ is a graph with two vertices
and three edges inside $T$,
having an open disc as a complement. That is, it is a trivalent spine of $T$.
Dually, this is a one-vertex triangulation of $T$.

\begin{figure}
\centering
\includegraphics[width = 8cm]{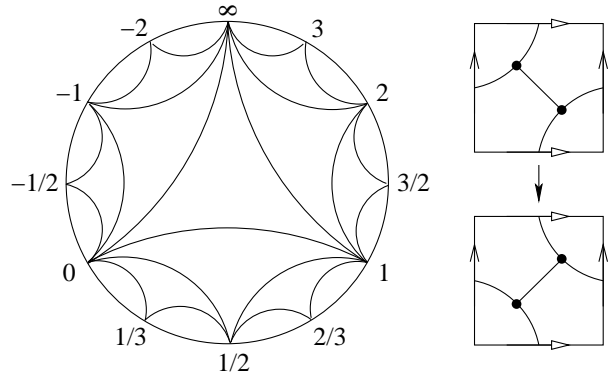} 
\caption{The Farey tesselation of the Poincar\'e disc into ideal triangles (left) and a flip (right).}
\label{tesselation:fig}
\end{figure}

The set of all $\theta$-graphs in $T$ up to isotopy can be described as follows. 
After choosing a meridian and a longitude,
every \emph{slope} on $T$ (\emph{i.e.}~isotopy class of simple closed
essential curves) is determined by a number $p/q\in\matQ\cup\{\infty\}$. 
Those numbers are the ideal vertices of the Farey tesselation
of the Poincar\'e disc sketched in Fig.~\ref{tesselation:fig}-left.
A $\theta$-graph contains three
slopes, which are the vertices of an ideal triangle of the tesselation. 
This gives a correspondence between the
$\theta$-graphs in $T$ and the triangles of
the tesselation.
Two $\theta$-graphs correspond to two adjacent triangles
when they share two slopes, \emph{i.e.}~when they are related by a \emph{flip}, 
shown in Fig.~\ref{tesselation:fig}-right.

\subsection{Manifolds with marked boundary}
Let $M$ be a connected compact 3-manifold with (possibly empty)
boundary consisting of tori. By associating to each torus component of
$\partial M$ a $\theta$-graph, we get a
{\em manifold with marked boundary}. 

Let $M$ and $M'$ be two marked manifolds, and $T\subset \partial M, T'\subset \partial M'$ 
be two boundary tori. 
A homeomorphism $\psi:T\to T'$ sending the marking of $T$ to the one of $T'$
is an \emph{assembling} of $M$ and $M'$. The result is a new marked manifold $N=M\cup_\psi M'$.
We define analogously a \emph{self-assembling} of $M$ along two tori $T,T'\subset\partial M$,
the only difference is that for some technical reason we allow the map to send one
$\theta\subset T$ either to $\theta'\subset T$ itself or to one of the $3$ other
$\theta$-graphs obtained from $\theta'$ via a flip. 

\subsection{Spines and complexity for marked manifolds}
The notion of spine extends from the class of closed manifold to the class 
of manifolds with marked boundary. 
\begin{defn}
Recall from Subsection~\ref{definitions:subsection} 
that a compact polyhedron is \emph{simple} when the link of each point is contained
in \includegraphics[width =.3 cm]{mercedes_small.eps}.
A sub-polyhedron $P$ of a manifold with marked boundary $M$ is
called a {\em spine} of $M$ if:
\begin{itemize}
\item $P\cup\partial M$ is simple,
\item $M\setminus(P \cup \partial M)$ is
an open ball,
\item $P \cap \partial M$ is a graph contained in the
marking of $\partial M$. 
\end{itemize}
\end{defn}
Note that $P$ is not in general a spine of $M$ in the usual 
sense\footnote{To avoid confusion, the term \emph{skeleton} was used in~\cite{MaPe}.}.
The {\em complexity} of a 3-manifold with marked boundary $M$ is of
course defined as the minimal number of vertices of a simple spine of 
$M$. Three fundamental properties extend from the closed case to the case with marked boundary:
complexity is still additive under connected sums, it is finite-to-one on
orientable irreducible manifolds,
and every orientable irreducible $M$ with $c(M)>0$ has a
minimal special spine~\cite{MaPe}.
(Here, a spine $P\subset M$ is \emph{special} when $P\cup\partial M$ is: the spine $P$ is
actually a \emph{special spine with boundary}, with $\partial P=\partial M\cap P$ 
consisting of all the
$\theta$-graphs in $\partial M$.)

\subsection{Bricks} 
An important easy fact is that if $M$ is obtained by assembling $M_1$ and $M_2$, and
$P_i$ is a spine of $M_i$, then $P_1\cup P_2$ is a spine of $M$. This implies
the first part of the following result.

\begin{proposition}[Martelli-Petronio~\cite{MaPe}] If $M$ is obtained by assembling
$M_1$ and $M_2$, we have $c(M)\leqslant c(M_1)+c(M_2)$. If $M$ is obtained by
self-assembling $N$, we have $c(M)\leqslant c(N)+6$.
\end{proposition}

When $c(M)=c(M_1)+c(M_2)$ or $c(M)=c(N)+6$, and the manifolds involved are 
irreducible\footnote{This hypothesis is actually determinant only in one case, see~\cite{MaPe}.},
the (self-)assembling is called \emph{sharp}.
\begin{defn} An orientable irreducible marked manifold $M$ is a \emph{brick} when
it is not the result of any sharp (self-)assembling.
\end{defn}

\begin{theorem}[Martelli-Petronio~\cite{MaPe}]
Every closed orientable irreducible $M$ is obtained from some bricks
via a combination of sharp (self-)assemblings.
\end{theorem}
There are only a finite number of such combinations giving the same $M$. 

\subsection{The algorithm that finds the bricks}
The algorithm described in Subsection~\ref{branches:subsection} also works 
for classifying all bricks
of increasing complexity, with some modifications, which we now sketch.
As we said above, every brick with $c>0$ has a minimal spine $P$ such that $P\cup\partial M$
is special. The $4$-valent graph $H=S(P\cup\partial M)$ contains the $\theta$-graphs marking
the boundary $\partial M$. By substituting (\emph{i.e.}~identifying) in $H$ 
each $\theta$-graph
with a point, we get a simpler $4$-valent graph $G$. We mark
the edges of $G$ containing that new points with a symbol $\star$.
It is then possible to encode the whole $P$ by assigning labels in $S_3$ on
the remaining edges of $G$, as in Subsection~\ref{branches:subsection}. The spine $P$
is uniquely determined by such data.

Every edge of $G$ can have a label in $S_3\cup\{\star\}$, giving $7^{2n}$ possibilities
to analyze during step (1b) of the algorithm
(actually, they are $2^n(3+1)^{2n}$ by Remark~\ref{ografi:rem}).
Although there are more possibilities to analyze than in the closed case
($7^{2n}$ against $6^{2n}$),
the non-minimality criteria for bricks listed below
are so powerful, that step (1b) is actually experimentally much quicker for bricks 
than for closed manifolds. This should be related with the experimental fact that
there are much more manifolds than bricks.

\begin{proposition}[Martelli-Petronio~\cite{MaPe}]
Let $P$ be a minimal special spine of a brick with $c>3$. 
The $3$ faces incident to an edge $e$ of $P$ are all distinct.
A face can be incident to at most one $\theta$-graph in $\partial P$.
\end{proposition}

\begin{figure}
\centering
\includegraphics[width = 8 cm]{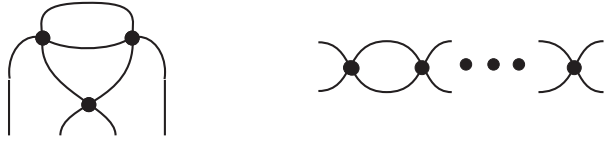} 
\caption{If $4$ edges disconnect $G$, then one of the two pieces is of one of these types.}
\label{2bridge:fig}
\end{figure}

\begin{theorem}[Martelli-Petronio~\cite{MaPe}] \label{bridge:teo}
Let $G$ be the $4$-valent graph associated to a minimal special spine  of a brick with $c>3$.
Then:
\begin{enumerate}
\item
no pair of edges disconnects $G$;
\item
if $c\leqslant 10$ and 
a quadruple of edges disconnects $G$, one of the two resulting components must be of one
of the forms shown in Fig.~\ref{2bridge:fig}.
\end{enumerate}
\end{theorem}
Point 2 of Theorem~\ref{bridge:teo} is proved for $c\leqslant 9$ in~\cite{MaPe} and conjectured
there to be true for all $c$: its extension to the case $c=10$ 
needed here is technical and we omit it.
We can restrict step (1b) of the algorithm to the \emph{useful}
$4$-valent graphs, \emph{i.e.}~the ones that are not forbidden by Theorem~\ref{bridge:teo}.
Table~\ref{useful:bricks:table} shows that only $2.1$ \% of the graphs are useful 
for $c=10,11$.

\begin{table}
\begin{center}
\begin{tabular}{cccccccccc} 
\phantom{\Big|} $n$    & $3$ & $4$  & $5$  & $6$  & $7$   & $8$    & $9$    & $10$    & $11$ \\
\hline
\phantom{\Big|} useful & $1$ & $2$  & $4$  & $11$ & $27$  & $57$   & $205$  & $1008$  & $6549$ \\
                all    & $4$ & $10$ & $28$ & $97$ & $359$ & $1635$ & $8296$ & $48432$ & $316520$
\end{tabular}
\end{center}
\caption{Useful graphs among all 4-valent graphs with $n\leqslant 11$ vertices.}
\label{useful:bricks:table}
\end{table} 

\subsection{Bricks with $c\leqslant 10$.} \label{bricks:found:subsection}
\begin{figure}[t]
\begin{minipage}{.03\textwidth}\hfil
\end{minipage}
\begin{minipage}{.2\textwidth}
\centering
\includegraphics[width = .9\textwidth]{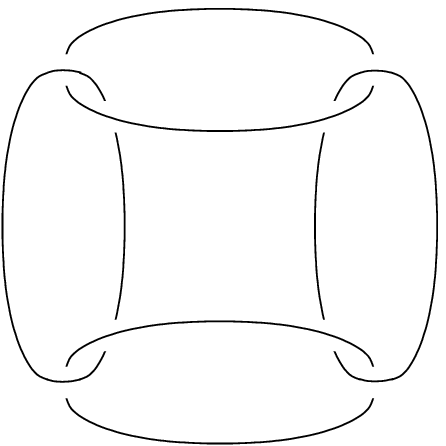}
\vglue -4mm
\end{minipage}
\begin{minipage}{.03\textwidth}\hfil
\end{minipage}
\begin{minipage}{.64\textwidth}
{The complement $M$ of this link
is a hyperbolic manifold. On each cusp,
there are two shortest loops of equal length,
and hence two preferred $\theta$-graphs, the ones containing
both loops. Up to symmetries of $M$,
there are only $3$ marked $M$'s
with such preferred $\theta$-graphs, and these are the ones with $c=10$.
}
\end{minipage}
\caption{The complement of a chain link with $4$ components.}
\label{chain4:fig}
\end{figure}

We list here the bricks found. There are two kinds of bricks:
the closed ones, and the ones with boundary. The closed ones correspond to the
closed irreducible 3-manifolds that do not decompose. 
\begin{theorem} \label{closed:bricks:teo}
The closed bricks having $c\leqslant 10$ are:
\begin{itemize}
\item
$\big(S^2,(2,1),(3,1),(m,1),-1\big)$ with $m\geqslant 5, m\neq 6$, having $c=m$;
\item
$\big(S^2,(2,1),(n,1),(m,1),-1\big)$ not of the type above and
with $\{n,m\}\neq\{3,6\},\{4,4\}$, having $c=n+m-2$;
\item
the closed hyperbolic manifold N.34 from Table~\ref{hyperbolic:table}, with volume $1.75712\ldots$
and homology $\matZ_7$, obtained as a $(1,-5,-3/2)$-surgery on the chain link, having $c=10$.
\end{itemize}
\end{theorem}
\begin{remark}
The manifolds $\big(S^2,(2,1),(n,1),(m,1),-1\big)$ with $\{n,m\}=\{3,6\}$ or $\{4,4\}$
are not bricks. Actually, they are flat torus bundles, whereas every other such manifold
is atoroidal.
\end{remark}
In the following statement, we denote by $N(\alpha,\beta,\gamma)$ the following marked manifold:
take the chain link of Fig.~\ref{chainlink:fig};
if $\alpha\in\matQ$, perform an $\alpha$-surgery on
one component, and if $\alpha = \theta^{(i)}$, drill that component and mark the new torus
with the $\theta$-graph containing the slopes $\infty, i$, and $i+1$. Do the same for $\beta$ and 
$\gamma$ (the choice of the components does not matter, see Fig.~\ref{chainlink:fig}).
\begin{theorem} \label{non:closed:bricks:teo}
The bricks with boundary having $c\leqslant 10$ are:
\begin{description}
\item[$c=0$:] one marked $T\times [0,1]$ and two marked solid tori;
\item[$c=1$:] one marked $T\times [0,1]$;
\item[$c=3$:] one marked (pair of pants)$\times S^1$;
\item[$c=8$:] one marked $\big(D,(2,1),(3,1)\big)$, and $N(1,-4,\theta^{(-1)})$;
\item[$c=9$:] four bricks of type $N(\alpha, \beta, \gamma)$, 
with $(\alpha,\beta,\gamma)$ being one of the following:
$$(1,-5, \theta^{(-1)}), \ (1, \theta^{(-2)}, \theta^{(-2)}), \
(\theta^{(-3)}, \theta^{(-2)}, \theta^{(-2)}), \ (\theta^{(-2)}, \theta^{(-2)}, \theta^{(-2)});$$
\item[$c=10$:] eleven bricks of type $N(\alpha, \beta, \gamma)$, 
with $(\alpha,\beta,\gamma)$ being one of the following:
$$(1,2,\theta^{(i)})\ {\rm with}\ i\in\{-3,-2,-1,0\}, \quad 
(1,-6, \theta^{(-1)}), $$ 
$$(-5, \theta^{(-2)}, \theta^{(-1)}), \quad 
(-5, \theta^{(-1)}, \theta^{(-1)}), \quad
(1, \theta^{(-1)}, \theta^{(-1)}), $$
$$(1, \theta^{(-4)}, \theta^{(-1)}), \quad
(2, \theta^{(-2)}, \theta^{(-2)}), \quad
(\theta^{(-3)}, \theta^{(-1)}, \theta^{(-1)}),$$
and three marked complements of the same link, shown in Fig.~\ref{chain4:fig}.
\end{description}
\end{theorem}

\begin{remark}
Using the bricks with $c\leqslant 1$, one constructs every marked solid torus. This
construction is the \emph{layered solid torus} decomposition~\cite{Bu, JaRu}.
An atoroidal manifold with $c\leqslant 10$ is either itself a brick,
or it decomposes into one brick $B$ of Theorem~\ref{non:closed:bricks:teo}
and some layered solid tori. 
\end{remark}

\begin{remark}
The generic graph manifold decomposes into some Seifert bricks with $c\leqslant 3$.
As Theorem~\ref{non:hyperbolic:teo} suggests, the only exceptions with $c\leqslant 10$
are the closed bricks listed by Theorem~\ref{closed:bricks:teo}, and some surgeries of the 
Seifert brick with $c=8$. 
\end{remark}
\begin{remark}
Table~\ref{hyperbolic:table} is deduced from Theorems~\ref{closed:bricks:teo}
and~\ref{non:closed:bricks:teo}, using SnapPea via a python script available from~\cite{weblist}.
\end{remark}
\begin{remark}
The proof of Theorem~\ref{non:hyperbolic:teo} from~\cite{MaPe:geometric} extends to $c=10$.
One has to check that the new hyperbolic bricks with $c=10$ do not contribute to
the complexity of non-hyperbolic manifolds, at least for $c=10$: we omit this discussion.
\end{remark}

We end this Section with a conjecture, motivated by our experimental results, which
implies that the decomposition into bricks is always finer than the \JSJ.
\begin{conjecture}
Every brick is atoroidal.
\end{conjecture}

\clearpage
\begin{finalspacing}{}

\bigskip

\name{Bruno Martelli}

\address{Dipartimento di Matematica \\%
Universit\`a di Pisa \\%
Via F.~Buonarroti 2 \\%
56127 Pisa, Italy%
}

\email{martelli@dm.unipi.it}

\bigskip
\classmark{57M27 (primary), 57M20, 57M50 (secondary)}


\label{martelli_last}

\end{finalspacing}

\end{document}